\documentclass[12pt]{article}
\usepackage{amssymb}
\usepackage{amsmath}
\usepackage{amsthm}
\usepackage{amsbsy}
\usepackage{graphicx}
\usepackage[]{epstopdf}
\newcommand{\be}{\begin{equation}}
\newcommand{\ee}{\end{equation}}
\newcommand{\bea}{\begin{eqnarray}}
\newcommand{\eea}{\end{eqnarray}}
\newcommand{\ba}{\begin{array}}
\newcommand{\ea}{\end{array}}

\newcommand{\bc}{\begin{center}}
\newcommand{\ec}{\end{center}}
\newcommand{\ben}{\begin{enumerate}}
\newcommand{\een}{\end{enumerate}}
\newcommand{\bfi}{\begin{figure}}
\newcommand{\efi}{\end{figure}}

\newcommand{\bq}{\begin{quote}}
\newcommand{\eq}{\end{quote}}
\newcommand{\bqu}{\begin{quotation}}
\newcommand{\equ}{\end{quotation}}
\newenvironment{emphit}{\begin{itemize}}{\end{itemize}}
\newcommand{\bemp}{\begin{emphit}}
\newcommand{\eemp}{\end{emphit}}

\newcommand{\bt}{\begin{tabular}}
\newcommand{\et}{\end{tabular}}

\newtheorem{myth}{Theorem}[section]
\newtheorem{mylem}{Lemma}[section]
\newtheorem{mycor}{Corollary}[section]
\newtheorem{mypro}{Proposition}[section]
\newtheorem{mydef}{Definition}[section]
\newtheorem{myrem}{Remark}[section]
\newtheorem{myexam}{Example}[section]

\begin{document}
\date{}
\title{The middle translations of finite involutory latin quandles
\footnote{2010 mathematics subject classification primary 20N05; secondary 57M27.}
\thanks{{\bf keywords: Involutory latin quandles, middle translations, representations, r-spin, l-spin,}}}
\author{A. O. Isere\thanks{All correspondence to be addressed to this author.}\\
Department of Mathematics,\\
Ambrose Alli University,\\
Ekpoma 310001, Nigeria.\\
isere.abed@gmail.com,\\
isereao@aauekpoma.edu.ng}\maketitle
\begin{abstract}
This paper studies the left (right) middle translations on finite involutory latin quandles and their representations. It also shows that if $(Q,\star)$ is a cyclic group of odd order n, such that
$x + y=L_{1}(x)\star \lambda_{1}(y)\star x$ for all $x,y\in Q$, where $L_{1}$ is a left translation, $\lambda_{1}$
a left middle translation and $1$ the identity element of the group $(Q,\star)$, then, $(Q,+)$ is a left involutory latin quandle of odd order n. Furthermore, the concept of spins of involutory latin quandle is investigated, and it is shown that if $Q(\cdot)$ is a left involutory latin quandle of odd order n, then the set of all right spins (r-spins) $\Phi_{R}$ is a cyclic group of odd order n under composition of mapping.
\end{abstract}
\section{Introduction}
\paragraph{}
The mappings $L_{a}$ and $R_{a}$ on a groupoid $(Q,\cdot)$ such that $L_{a}:Q\rightarrow Q$ and $R_{a}:Q\rightarrow Q$ ($a\in Q$) defined as $L_{a}(x)=a\cdot x$ and $R_{a}(x)=x\cdot a$ for all $x\in Q$ are called left and right translations respectively. These translations find relevant applications in quasigroups. Quickly, a groupoid $(Q,\cdot)$ is called a quasigroup if the translations $L_{a}$ and $R_{a}$ are bijective. For a finite quasigroup $Q$ the translations $L_{a}$ and $R_{a}$ are permutations. On the other hand, the permutations $\lambda_{i}, \varphi_{i} (i\in Q)$ of $Q$ that are defined as
$$
\lambda_{i}(x)\cdot x=i
$$
and
$$
x\cdot \varphi_{i}(x)=i
$$
for all $x \in Q$ are called left and right middle translations of an element $i$ in a quasigroup $Q(\cdot)$ respectively (see-\cite{qua15, qua28}). Moreover, the translation that is both a left and right middle translation is simply called the middle translation. These translations were first introduced and studied by Belousov (\cite{qua11,qua12}), and since then
many researchers have developed interest in expanding the concept (see-\cite{qua13,qua15,qua14,qua28,qua30}). Therefore, the focus of this paper is to investigate the middle translations and their representations in relation to finite involutory latin quandles. The term involutory latin quandles, as used in this paper, refers to latin quandles with involutory properties. Detail study of quandles and involutory latin quandles abound in  literature(see-\cite{qua19,qua21,qua02,qua27a,qua27}).
\par
The permutation $\lambda_{i}(x)$ means left track (l-track). That is, to find in the column $x$ the cell containing an element $i$, we must select the row $\lambda_{i}(x)$ \cite{qua14}. In other words, $\lambda_{i}(x)$ is a permutation
of those elements that multiply $x$ from the left to give the result $i$. Thus $\lambda_{i}(x)$ is a row selection. On the other hand, $\varphi_{i}(x)$ is a column selection, and a permutation of those elements that multiply $x$ from the right
to give $i$. Therefore $\varphi_{i}(x)$ is a right track (r-track). The group generated by all these
translations ($R_{i}, L_{i}, \lambda_{i} ~\&~ \varphi_{i}$) is called the multiplication group denoted as $M(Q,\cdot)$
where $(Q,\cdot)$ is a quasigroup.
\begin{mydef} (\cite{qua01,qua09})
A quandle is a set $X$ with a binary operation $(a,b)\mapsto a b$ such that
\begin{description}
\item[(1)] For any $x\in X,~ xx=x$
\item[(2)] For any $a,b\in X$, there is a unique $x\in X$ such that $a=xb$
\item[(3)] For any $a,b,c\in X$, $(a b) c = (a c)(b c)$
\end{description}
\end{mydef}
 The juxtaposition represents the binary operation in most of the definitions in this paper.
\begin{mydef} (\cite{qua16,qua23})
A quandle $X$ is commutative if it satisfies the identity
$$ x y=y x ~~ \forall ~ x,y\in X.$$
\end{mydef}
\begin{mydef} (\cite{qua18,qua01})
An abelian quandle is a quandle satisfying the identity:
$$(w x) (y z)=(wy) (x z)$$
\end{mydef}
\begin{mydef} (\cite{qua17})
Given two quandles $(X,\star)$ and $(Y,\bullet)$ , a map $f:(X,\star)\rightarrow (Y,\bullet)$
is a quandle homomorphism if
$$
f(a\star b)=f(a)\bullet f(b)~~\forall~a,b\in X
$$
If $f$ is a bijection then $f$ is called an isomorphism, and  $(X,\star)$ and $(Y,\bullet)$ are said to be isomorphic quandles.
\end{mydef}
\begin{mydef}
The automorphism group of quandle $(X, *)$ , denoted as $Aut(X)$ is the group of all
isomorphisms $f: X \rightarrow X$.
\end{mydef}
\begin{mydef} (\cite{qua16,qua07})
The inner automorphism group of a quandle $(X,\star)$ denoted as Inn$(X)$ is the subgroup of Aut$(X)$ generated by all $S_x$, where $S_{x}(y)=y\star x$, for any $x,y\in X$. The map $S_{x}: X\rightarrow X$ that maps $u$ to $u\star x$ defines a right action of $X$ on $X$, so that we obtain a map $X\rightarrow Inn(X)$
\end{mydef}
\begin{mydef} \cite{qua11,qua19}\label{CoreDef}

Let $(G,\cdot)$ be a group, or, more generally, a Bol loop. The binary algebra $(G,\star)$ with
 \begin{align}\label{quad16b}
 x\star y=xy^{-1}x
 \end{align}
 is an involutory quandle, called the core of $(G,\cdot)$.
\end{mydef}

 Bruck \cite{phd41} had earlier shown that the core of a Moufang loop, originally defined as
\begin{align}\label{quad16}
x+y=yx^{-1}y
\end{align}
is an involutory quandle.

\par It is to be noted that the quandles described above are not quasigroups or latin quandles.
\begin{mydef}\cite{qua19}
A groupoid $(Q,\star)$ is called a latin quandle if it obeys the following laws simultaneously:
\begin{description}
\item[(i)] $x\star x=x$ for all $x\in Q$ [idempotent law]
\item[(ii)] $a\star x=b$ for all $x\in Q$ and $a,b$ specified in $Q$ [left division law]
\item[(iii)]$y\star a=b$ for all $y\in Q$ and $a,b$ specified in $Q$ [right division law]
\item[(iv)] $ a\star (x\star y)=(a\star x)\star (a\star y)$ for all $a, x$ and $y$ in $Q$ [left distributive law]
\item[(v)] $(x\star y)\star a =(x\star a)\star (y\star a)$ for all $a, x$ and $y$ in $Q$ [right distributive law]
\end{description}
\end{mydef}
\begin{myexam}\label{exam1}
Let $(G,\cdot)$ be a cyclic group (more generally, a commutative Moufang loop) of odd order $n$ such that
$$
x+y=xy^{-1}\cdot x ~~\forall~x,y\in G.
$$
Then, the core $(G,+)$ is a latin quandle of odd order $n$
\end{myexam}
\begin{mydef}(\cite{qua21})\label{quad19}
A latin quandle $(Q,\circ)$ that obeys properties:
\begin{description}
\item[(1)] $x\circ(x\circ y)=y$ Left Involutory Property (LIP) is called a LIPQ.
\item[(2)] $(y\circ x)\circ x =y$ Right Involutory Property (RIP) is called a RIPQ.
\item[(3)](1) and (2)  Involutory Property (IP) is called an IPQ.
\item[(4)] $ x\circ (y\circ x) =y$ or $y=(x\circ y) \circ x$ Cross Involutory Property (CIP) is called a CIPQ.
\end{description}
For all $x,y\in Q$.
\end{mydef}
Example \ref{exam1} above is a LIPQ of odd order $n$.
In a similar consideration, (\ref{quad16}) also gives the core $(G,+)$ as a RIPQ of odd order $n$.
Therefore, involutory latin quandles can be constructed as cores of cyclic groups as presented above.
\par
A latin quandle $Q$ that obeys Definition \ref{quad19} (LIPQ, RIPQ or IPQ) is the focus of this paper.
These acronyms henceforth will be used to represent these involutory quandles, and are not new in literature. They have always been used in connection with inverse properties quasigroups (loops)as shown in the following definition.
\begin{mydef}\label{quad18} (\cite{qua28})
\begin{enumerate}
\item
A quasigroup $(Q,\circ)$ has the Left Inverse Property (LIP) if there exists a permutation $\lambda$ of the set $Q$ such that $$ \lambda x \circ (x\circ y)=y$$ for all $x,y\in Q$ (By Belousov).
\item
A quasigroup $(Q,\circ)$ has the Right Inverse Property (RIP) if there exists a permutation $\rho$ of the set $Q$ such that $$   (x\circ y)\circ \rho y= x$$ for all $x,y\in Q$ (By Belousov).
\item
A quasigroup $(Q,\circ)$ has the Inverse Property (IP) if it is  a LIP and RIP-quasigroup (By Belousov).
\item
A quasigroup $(Q,\circ)$ has the Cross Inverse Property (CIP) if there exists a permutation $J$ of the set $Q$ such that $$   (x\circ y)\circ Jx= y$$ for all $x,y\in Q$ (First by Artzy and later by Keedwell and Shcherbacov) (\cite{qua22,qua24,qua28}).
\end{enumerate}
\end{mydef}
More results on involutory latin quandles (LIPQ, RIPQ, CIPQ, IPQ) can be found in \cite{qua21}. A few are presented below.
\begin{myth}\cite{qua21} Let $Q(\cdot)$ be a cyclic group of odd order n such that
$$ x+y=(y^{-1}x)x, \forall x,y\in Q.$$ Then, $(Q,+)$ is a LIPQ of odd order n.
\end{myth}
\begin{myth}\cite{qua21} Let $Q(\cdot)$ be a cyclic group of odd order n such that
$$ x+y=y(yx^{-1}), \forall x,y\in Q.$$ Then, $(Q,+)$ is a RIPQ of odd order n.
\end{myth}
\begin{myth}\cite{qua21} Let $Q(\cdot)$ be a commutative group (not cyclic) of order $3^{n}, n\ge 1$ such that
$$ x+y=x(y^{-1}x), \forall x,y\in Q.$$ Then, $(Q,+)$ is an IPQ of order $3^{n}, n\ge 1$.
\end{myth}
The authors also stated that a latin Alexander quandle of order $4^{n}, n\ge 1$ is a CIPQ, and presented several concrete non-isomorphic examples of LIPQs, RIPQs, CIPQs and IPQs (see-\cite{qua21}).

\begin{mydef} (\cite{phd20})
Let $Q$ be a quasigroup (loop). The set $\Pi=\{ R(a):a\in Q \}$ is called the right regular representation
of $Q$ or briefly the representation of $Q$. The left regular representation is defined analogously.
\end{mydef}
\begin{mydef} (\cite{qua14})
By a spin of a quasigroup $Q(\cdot)$ we mean the permutation $$\varphi_{ij}=\varphi_{i}\varphi^{-1}_{j}=\varphi_{i}\lambda_{j}$$ where $\varphi_{i}$ and $\lambda_{j}$ are
tracks of $Q(\cdot)$. The spin $\varphi_{ii}$ is called trivial. And the set of all spins of a quasigroup $Q(\cdot)$
is denoted $\Phi_{Q}(\cdot)$.
\end{mydef}
\begin{mycor}\label{quad02a}\cite{qua14}
 In any group $G(\cdot)$ we have
 \begin{enumerate}
 \item $\varphi_{i}(x)=x^{-1}\cdot i  \quad (\lambda_{i}(x)=i\cdot x^{-1})$,
 \item $\varphi_{1}(x)=\lambda_{1}(x)=x^{-1}$
 \item $L_{i}(x)=\lambda_{i}(x)\cdot x^{2}$
 \item $R_{i}(x)=x^{2}\cdot \varphi_{i}(x)$
 \end{enumerate}
 where $1$ is the identity element of the group $G$.
 \end{mycor}
 The next section presents the left and right middle translations and their induced representations as well as the construction of involutory latin quandles using the middle translations (see Theorem \ref{quad10}  and Theorem \ref{quad11} ). Section 3 presents the applications of the middle translations in l-spins and r-spins. These concepts enable us to recover cyclic groups (see Theorem \ref{quad04}  and Theorem \ref{quad12}) from involutory latin quandles.

\section{Middle translations and their representations}
\begin{mylem}\label{quad19}
Let $Q (\cdot)$ be an involutory latin quandle (LIPQ, RIPQ or CIPQ), then the following hold:
\begin{description}
\item[(1)] $\lambda_{i}=\varphi^{-1}_{i}$
\item[(2)] $\varphi^{-1}_{i}(x)\cdot x =i$
\item[(3)] $L_{i}(x)=(\lambda_{i}(x)\cdot x)\cdot x$
\item[(4)] $ L_{i}(x)=( x\cdot \varphi_{i}(x))\cdot x$
\item[(5)] $R_{i}(x)= x\cdot (\lambda_{i}(x)\cdot x)$
\item[(6)] $R_{i}(x)=x\cdot (x\cdot \varphi_{i}(x))$
\end{description}
For all $x\in Q$ and $i\in Q$.
\end{mylem}
The Proof of the above results simply follow from the definitions of the four translations on an involutory latin quandle $Q$.
\begin{mylem}\label{quad01}
Let $Q(+)$ be a finite RIPQ. Then, for each $i\in Q, \lambda_{i}(x)+x=i$ implies that $i+x=\lambda_{i}(x)$  for all $ x\in Q$.
\end{mylem}
{\bf Proof:}\\
Consider a RIPQ, then $(i+x)+x=i$ for all $i,x\in Q$. Also $$\lambda_{i}(x)+x=i \Rightarrow (i+x)+x=\lambda_{i}(x)+x=i$$ Thus $$ i+x=\lambda_{i}(x)$$
\begin{mylem}\label{quad02}
Let $Q(\cdot)$ be a finite LIPQ. Then, for each $i\in Q, x\cdot \varphi_{i}(x)=i$ implies that $x\cdot i=\varphi_{i}(x)$  for all $ x\in Q$.
\end{mylem}
{\bf Proof:}\\
$LIPQ$ guarantees that
 $$x\cdot (x\cdot i)=i$$ and by definition $$  x\cdot \varphi_{i}(x)=i$$ Thus $$ x\cdot i=\varphi_{i}(x)$$

\begin{mycor}\label{quad08}
Let $Q(\cdot)$ be a latin quandle. Then, the permutations $L_{i}$ and $\lambda_{i}$ coincide if and only if $Q(\cdot)$ is a RIPQ.
\end{mycor}
{\bf Proof:}\\
Given that $Q(\cdot)$ is a RIPQ and that $L_{i}:x\mapsto i\cdot x$ ($i\in Q$), then by Lemma \ref{quad01} the first part holds.\\
Conversely, suppose that $L_{i}=\lambda_{i}$, then $i\cdot x=\lambda_{i}(x)$. Multiplying both sides by $x$ from the right gives
$(i\cdot x)\cdot x=\lambda_{i}(x)\cdot x=i  ~ ~ \forall  x\in Q$. Thus, $Q(\cdot)$ is a RIPQ.

\begin{mycor}\label{quad09}
Let $Q(\cdot)$ be a latin quandle. Then, the permutations $R_{i}$ and $\varphi_{i}$ on $Q(\cdot)$ defined as $R_{i}(x)=x\cdot i$ and $x\cdot \varphi_{i}(x)=i$ respectively coincide if and only if $Q(\cdot)$ is a LIPQ.
\end{mycor}
{\bf Proof:}\\
 Since $Q(\cdot)$ is a LIPQ and that $R_{i}:x\mapsto x\cdot i$ ($i\in Q$), then by Lemma \ref{quad02} the forward statement holds.\\
The converse is proved as above.
\begin{mypro}\label{quad05}
Let $Q$ be a finite set and assume there exists a map $\lambda_{i}:Q\rightarrow Q $ ($i\in Q$). Then the binary operation $\circ$ defined by $\lambda_{i}(x)\circ x=i$ gives a RIPQ if and only if:
\begin{enumerate}
\item
$\lambda_{x}(x)=x$ for all $x\in Q$
\item
$\lambda_{i}$ is bijective
\item
$\lambda_{(i\circ x)}(y)=\lambda_{i}(y)\circ \lambda_{x}(y)$ for all $x,y \in Q$.
\end{enumerate}
\end{mypro}
{\bf Proof:}\\
Given that $(Q,\circ)$ is a RIPQ, then
\begin{enumerate}
\item $x\circ x=x$ implies that $\lambda_{x}(x)=x $ by Lemma \ref{quad01}.
\item By Corollary \ref{quad08}, the permutation $\lambda_{i} $ is bijective since $Q$ is a latin quandle.
\item  Lemma \ref{quad01} and distributivity guarantees $\lambda_{(i\circ x)}(y)=(i\circ x)\circ y =\lambda_{i}(y)\circ \lambda_{x}(y)$.
\end{enumerate}
The converse follows from a careful application of Lemma \ref{quad01} and Corollary \ref{quad08}
\begin{mypro}\label{quad05b}
Let $Q$ be a finite set and assume there exists a map $\varphi_{i}: Q\rightarrow Q$ ($i\in Q$). Then the binary operation $\circ$ defined by $ x\circ \varphi_{i}(x)=i $ gives a LIPQ if and only if:
\begin{enumerate}
\item
$\varphi_{x}(x)=x $ for all $x\in Q$
\item
$\varphi_{i}$ is bijective
\item
$\varphi_{(i\circ x)}(y)=\varphi_{i}(y)\circ \varphi_{x}(y)$ for all $x,y \in Q$.
\end{enumerate}
\end{mypro}
{\bf Proof:}\\
The first part of the proof is similar to the proof of Proposition \ref{quad05} and the converse also follows from Lemma \ref{quad02} and Corollary \ref{quad09}
\begin{mycor}
Let $Q$ be a RIPQ such that $\lambda_{i}:Q\rightarrow Q $ ($i\in Q$). Then $\lambda_{i}$  defined as $\lambda_{i}(x)\circ x=i$ is an automorphism.
\end{mycor}
{\bf Proof:}\\
The definition $\lambda_{i}(x)\circ x=i$ implies that $\lambda_{i}(x)=i\circ x$ by Lemma \ref{quad01}.\\
Consider: $\lambda_{i}(x\circ y)=i\circ (x\circ y)=\lambda_{i}(x)\circ \lambda_{i}(y)$.
Thus, $\lambda_{i}$ is an automorphism since $Q$ is bijective as a latin quandle.
\begin{myrem}
That $\varphi_{i}$ is an automorphism when $(Q)$ is a LIPQ can be proved in a similar manner.
\end{myrem}

\begin{mydef}\label{quad06}
Let $Q (\cdot)$ be an involutory latin quandle. The sets $\Pi_{\lambda}=\{ \lambda_{i}(x): \lambda_{i}(x)\cdot x=i, x\in Q \}$ and $\Pi_{\varphi}=\{ \varphi_{i}(x): x \cdot \varphi_{i}(x)=i, x\in Q \}$ are called the left and right middle representations respectively.
\end{mydef}
\begin{mypro}\label{quad07}
A set $\Pi$ of permutations on $Q$ is the representation of left middle translations on an involutory latin quandle $Q(\cdot)$ if and only if
\begin{description}
\item[(i)] $\pi_{x}(x)=x$ for all $x\in Q$ and $\pi_{x} \in \Pi$ (i.e every permutation fixes an element of $Q$),
\item[(ii)] for all $x,y\in Q$, there exists a unique $\pi_{y} \in \Pi$ such that $\pi_{y}(x)\cdot x=y$,
\item[(iii)] $\alpha,\beta \in \Pi$ and $\alpha\beta$ fixes the same element of $Q$, then $\alpha=\beta$.
\end{description}
\end{mypro}
{\bf Proof:}\\
Suppose first that (i), (ii) and (iii) hold. Then we need to show that $\Pi$ is a set of permutations induced by the left
middle translations on a latin quandle $Q$.\\
That $\pi_{x}(x)=x$ and $\pi_{y}(x)\cdot x=y$ means that $\pi_{i}(x)\cdot x=i ~ \forall ~ x\in Q$. Moreover, if $\alpha$ and $\beta$ are in $\Pi$
and their product fixes the same element say $x$, then $\alpha \beta=\pi_{x}(x)=x \Rightarrow \pi_{x}(x)\cdot x=x$. Therefore, $\pi_{x} \in \Pi$ (by Definition \ref{quad06}). Thus, $\Pi$ is a set of permutations induced by the left middle translations.\\
Conversely,
\begin{description}
\item[(i)] given $\Pi$ as in Definition \ref{quad06}, $\pi_{i}(x)\cdot x=i (i\in Q)$. So, $i=x$ means that $\pi_{x}(x)\cdot x=x\cdot x$ (since $Q$ is a latin quandle). Then, $\pi_{x}(x)=x ~ \forall ~ x\in Q$.
\item[(ii)] If $i=y, \pi_{y}(x)\cdot x=y~ ~ \forall x, y\in Q $. Then to show uniqueness: suppose $\pi_{y}$ is not unique in $\Pi$, then there exists $\pi_{y}$ and $\pi_{y}'$ such that $\pi_{y}(x)\cdot x=y$ and $\pi_{y}'(x)\cdot x=y$. This implies that $\pi_{y}(x)\cdot x=\pi_{y}'(x)\cdot x \Rightarrow \pi_{y}=\pi_{y}'$. Thus $\pi_{y}$ is unique in $\Pi$.
\item[(iii)] Let $\alpha=\pi_{i}$ and $\beta=\pi_{x} ~(i\ne x)$. Then $ \alpha\beta(y)=\pi_{i}(y)\cdot \pi_{x}(y)=\pi_{(i\cdot x)}(y)$(Proposition \ref{quad05}).
Thus $\alpha\beta=\pi_{(i\cdot x)}$. From (i), $\pi_{(i\cdot x)}(i\cdot x)=(i\cdot x)$. So, $\alpha\beta$ fixes $(i\cdot x)~~ i\ne x$. But if $i=x$ then $\pi_{i}=\pi_{x}$ and $\alpha\beta$ fixes $x$, the same element as $\alpha$ and $\beta$. Thence, $\alpha=\beta$.
\end{description}
\begin{mypro}
A set $\Pi$ of permutations on $Q$ is the representation of right middle translations on an involutory latin quandle $Q(\cdot)$ if and only if
\begin{description}
\item[(i)] $\pi_{x}(x)=x$ for all $\pi_{x} \in \Pi$ (i.e every permutation fixes an element),
\item[(ii)] for all $x,y\in Q$, there exists a unique $\pi_{y} \in \Pi$ such that $x\pi_{y}\cdot x=y$,
\item[(iii)] $\alpha,\beta \in \Pi$ and $\alpha\beta$ fixes the same element of $Q$, then $\alpha=\beta$.
\end{description}
\end{mypro}
The Proof is similar to the proof of Proposition \ref{quad07}.
\begin{myth}
Let $Q(+)$ be a latin quandle of odd order n. Then, the representation induced by the  left middle translations on $Q(+)$ is a LIPQ of odd order n if and only if $Q(+)$ is a RIPQ  of odd order n.
\end{myth}
{\bf Proof:}\\
Let $Q(+)$ be a latin quandle generated by $(x+y)+y=x$ for all $x,y\in Q$. Then for each $x$ in $Q$ this implies  $\lambda_{x}(y)+ y=x$. This gives the left middle translation.\\
Conversely, suppose $Q(+)$ is generated by the set of all left middle translations on $Q(+)$. That is $$ \lambda_{i}(x)+ x=i $$ By Lemma \ref{quad01} $$ (i+ x)+x =i $$ This means that $Q(+)$ is a RIPQ.
\begin{myth}\label{quad03}
Let $Q(+)$ be a RIPQ of order $n$. Then, the representation induced by the left middle translation on $Q(+)$ is a CIPQ of order $n$ if and only if the representation is commutative.
\end{myth}
{\bf Proof:}\\
Suppose that the representation induced by the left middle translations on $Q(+)$ is commutative. Then $$ \lambda_{i}(x)+ x=i \Rightarrow x+ \lambda_{i}(x)=i$$ But $$ \lambda_{i}(x)=(i+ x)\Rightarrow x+ (i+ x)=i $$ This gives a CIPQ.\\
Conversely, given that $Q(+)$ is a CIPQ, that is, $x+(y+ x)=y$. Represent $y + x$ by $\lambda_{y}(x)$ then
$x+ \lambda_{y}(x)=x$. Commutativity implies that $\lambda_{y}(x)+ x=y$. This gives the left middle translation.

\begin{myth}\label{quad10}
Let $(Q, \star)$ be a cyclic group of odd order n such that $x + y=L_{1}(x)\star \lambda_{1}(y)\star x ~~\forall x,y\in Q$, where $L_{1}$ is a left translation, $\lambda_{1}$ a left middle translation and $1$ the identity element of $(Q,\star)$. Then, $(Q,+)$ is a LIPQ of  odd order n.
\end{myth}
{\bf Proof:}\\
 consider:
$$
(x + y)= L_{1}(x)\star \lambda_{1}(y)\star x
$$
then simplifying by, $ \lambda_{1}(y)=y^{-1}$ where $1$ is the identity element of $(Q,\star)$ ensures the following results: $$ \begin{gathered}
x + x= L_{1}(x)\star \lambda_{1}(x)\star x=x\\
(x + y)+ z=L_{1}[L_{1}(x)\star \lambda_{1}(y)\star x]\star \lambda_{1}(z)\star L_{1}(x)\star \lambda_{1}(y)\star x\\
=(x + z)+ (y + z)
\end{gathered}
$$
Similarly, $$ x+(y+z)=L_{1}(x)\star \lambda_{1}[L_{1}(y)\star \lambda_{1}(z)\star y]\star x=(x+y)+(x+z)$$
Left and  right divisibility hold since $(Q,\star)$ is a group. Then, $(Q,+)$ is a latin quandle.\\ Next
$$
\begin{gathered}
x+(x + y)= x + [L_{1}(x)\star \lambda_{1}(y)\star x]\\
=L_{1}(x)\star \lambda_{1}[L_{1}(x)\star \lambda_{1}(y)\star x]\star x=y
\end{gathered}
$$
Thus, $x + (x + y)=y$. Therefore, $(Q,+)$ is a LIPQ.
\begin{myth}\label{quad11}
Let $(Q,\star)$ be a cyclic group of odd order n such that $x + y=R_{1}(y)\star \varphi_{1}(x)\star y$ ($x,y\in Q$) where $R_{1}$ is a right translation, $\varphi_{1}$ a right  middle translation and $1$ is the identity element of $(Q,+)$. Then $(Q,+)$ is a RIPQ of odd order n.
\end{myth}
{\bf Proof:}\\
Similarly, for any group $(Q,\star)$, $\varphi_{1}(x)=x^{-1}$ where $1$ is the identity element of $(Q,\star)$, one can show that
$ x + x=x $ and $(x + y) + z = (x + z) + (y + z)$ as above. The left and right divisibility hold for the same reason. Therefore, $(Q,\star)$ is a latin quandle.\\ Next
$$
\begin{gathered}
(x + y) + y\\
=[R_{1}(y)\star\varphi_{1}{x}\star y] + y\\
=R_{1}(y)\star\varphi_{1}[R_{1}(y)\star \varphi_{1}{x}\star y]\star y=x
\end{gathered}
$$
Thus, $(x + y) + y=x$. Therefore, $(Q,+)$ is a RIPQ.
\begin{myth}\label{quad10a}
Let $(Q, \circ)$ be a commutative group (not necessarily cyclic) of order $3^{n}, n\ge 1$  such that $$x + y=L_{1}(x)\circ \lambda_{1}(y)\circ x ~~\forall x,y\in Q, $$ where $L_{1}$ is a left translation, $\lambda_{1}$ a left middle translation and $1$ the identity element of $(Q,\circ)$. Then, $(Q,+)$ is an IPQ of order $3^{n}, n\ge 1$.
\end{myth}
{\bf Proof:}\\
 consider:
$$
(x + y)= x \circ  \lambda_{1}(y)\circ L_{1}(x)
$$
then simplify by $ \lambda_{1}(y)=y^{-1}$ gives:
$$ \begin{gathered}
x + x= x \circ  \lambda_{1}(x)\circ L_{1}(x)=x\\
(x + y)+ z=[x \circ \lambda_{1}(y)\circ L_{1}(x)]+z\\ = (x \circ \lambda_{1}(y)\circ L_{1}(x))\circ \lambda_{1}(z)\circ L_{1}(x \circ \lambda_{1}(y)\circ L_{1}(x))\\= (xy^{-1}x)z^{-1}(xy^{-1}x)\\
=(x + z)+ (y + z)
\end{gathered}
$$
Similarly, $$ \begin{gathered}
x+(y+z)= x + [y \circ \lambda_{1}(z)\circ L_{1}(y)]\\
 = x\circ \lambda_{1}[y \circ \lambda_{1}(z)\circ L_{1}(y)]\circ L_{1}(x)\\ =(x+y)+(x+z)
 \end{gathered}
 $$
Then left and  right divisibility hold since $(Q,\circ)$ is a group. Therefore, $(Q,+)$ is a latin quandle.\\ Next
$$
\begin{gathered}
x+(x + y)= x + [x \circ \lambda_{1}(y)\circ L_{1}(x)]\\
= x\circ  \lambda_{1}[x \circ \lambda_{1}(y)\circ L_{1}(x)]\\ = \lambda_{1}[x \circ \lambda_{1}(y)\circ L_{1}(x)]\circ x= (y+x)+x
\end{gathered}
$$
Thus, $(Q,+)$ is an IPQ.
\begin{myrem}
$(Q,+)$ can equivalently be defined as $x + y=x\circ \varphi_{1}(y)\circ R_{1}(x) ~~\forall x,y\in Q, $ where $R_{1}$ is a right translation, $\varphi_{1}$ a right middle translation and $1$ the identity element of $(Q,\circ). $
\end{myrem}
\begin{myth}\label{quad21}
Let $Q(\circ)$ be a commutative latin quandle of order $3^{n}, n\ge 1$ then $Q$ is an IPQ of order $3^{n}$ if and only if
$\lambda_{i}= \varphi_{i}$.
\end{myth}
{\bf Proof}\\
Suppose $Q(\circ)$ is an IPQ, then $Q(\circ)$ is a LIPQ and RIPQ simultaneously. A nice application of Lemma \ref{quad01} and Lemma \ref{quad02} and since $Q$ is commutative ensures that
$$ (i\circ x)\circ x = x\circ (x\circ i)=i=\lambda_{i}(x)\circ x= \varphi_{i}(x)\circ x$$
implies that $\lambda_{i}=\varphi_{i}$.\\
Conversely, since $\lambda_{i}(x)=\varphi_{i}(x)$, then $\lambda_{i}(x)\circ x =\varphi_{i}(x)\circ x$
But $\lambda_{i}(x)\circ x = i=x\circ \varphi_{i}(x)$. Then $(i\circ x)\circ x = i=x\circ (x\circ i)$ (by Lemma \ref{quad01} and Lemma \ref{quad02}),  Thus $Q$ is an IPQ.
\begin{myrem}
All translations $L_{i}, R_{i}, \lambda_{i}$ and $\varphi_{i}$ coincide in IPQ.
\end{myrem}

\section{Spins of Involutory Latin Quandles}
\begin{mydef}
Let $Q(\cdot)$ be a latin quandle. Then, by a left spin (l-spin) of $Q(\cdot)$ we mean the permutation
$$ \lambda_{ij}=\lambda_{i}\lambda^{-1}_{j}=\lambda_{i}\varphi_{j} $$
where $\lambda_{i}$ and $\varphi_{j}$ are left and right middle translations on $Q$ respectively.
\end{mydef}
\begin{mydef}
Let $Q(\cdot)$ be a latin quandle. Then, by right spin (r-spin) of $Q(\cdot)$ we mean the permutation
$$ \varphi_{ij}=\varphi_{i}\varphi^{-1}_{j}=\varphi_{i}\lambda_{j} $$
where  $\varphi_{i}$ and $\lambda_{j}$ are right and left middle translations on $Q$ respectively.
\end{mydef}
A permutation $\pi_{ij}$ on a latin quandle $Q$ is a spin if it is both l-spin and r-spin.
\begin{mylem}\label{quad20}
Let $Q(\cdot)$ be an involutory latin quandle (LIPQ, RIPQ, or both) of order $n$. Then, the following properties hold:
\begin{enumerate}
\item $\varphi_{ij}(x)\neq x ~(\lambda_{ij}(x)\neq x)$ for all $x\in Q$ and $i\neq j$
\item $\varphi_{pi}(x) \neq \varphi_{pj}(x) (\lambda_{pi}(x) \neq \lambda_{pj}(x))$ for all $x\in Q$ and $i\neq j$
\item $\varphi_{ij}=\varphi^{-1}_{ji} (\lambda_{ij}=\lambda^{-1}_{ij})$ for $i\neq j$
\item $\varphi_{ij} = \varphi_{(i+1)(j+1)}$ for $i=1,2,...,n-1;j=1,2,...,n-1$ and $\varphi_{ii}$ is trivial.
\item $\lambda_{ij} = \lambda_{(i+1)(j+1)}$ for $i=1,2,...,n-1;j=1,2,...,n-1$ and $\lambda_{ii}$ is trivial.
\item $ \varphi_{n1}=\varphi_{(n-1)n} $
\item $\lambda_{n1}=\lambda_{(n-1)n}$
\end{enumerate}
\end{mylem}
\begin{myth}\label{quad04}
Let $Q(\cdot)$ be a LIPQ of odd order n. Then, the set of all r-spins of $Q$ is a cyclic group of odd order n under composition of mapping, denoted as $(\Phi_{R},\circ)$.
\end{myth}
{\bf Proof:}\\
Let $$ P_{R}=\{ \varphi_{i},\lambda_{j}\in M(Q,\cdot)|\varphi_{ij}=\varphi_{i}\varphi^{-1}_{j}=\varphi_{i}\lambda_{j},i,j\in Q \} $$ such that the order of $P_{R}$ is odd. From the definition, $P_{R}$ is a subset of the multiplication group $M(Q,\cdot)$. $\varphi_{ii}\in P_{R}$ since $\varphi_{ii}=\varphi_{i}\varphi^{-1}_{i}=I_{ii}$ (identity r-spin), and thus, $P_{R}$ is not empty. Also, $\varphi_{ij}\circ \varphi_{jk}=\varphi_{ik} (i\ne k)$. Hence $P_{R}$ is closed under composition of mapping. Moreover, $\varphi_{ii}=\varphi_{ij}\circ \varphi^{-1}_{ij}=\varphi_{ij}\circ \varphi_{ji}=I_{ii}$. Thus, $\varphi_{ji}$ is the inverse of
$\varphi_{ij}$ in $P_{R}$. Therefore, $P_{R}$ is a subgroup of $M(Q,\cdot)$, and thus a group.
\par
Now, consider
$\varphi_{ij}(x)=x\cdot (ij)=(x\cdot i)j=j(i\cdot x)=(ji)\cdot x=\varphi^{-1}_{ji}(x) \Rightarrow\varphi_{ij}=\varphi^{-1}_{ji}$ ( by Lemma \ref{quad20}(3)). Thus $P_{R}$ is commutative. Therefore, $P_{R}$ is a cyclic group of odd order n.
\begin{myth}\label{quad12}
Let $Q(\cdot)$ be a RIPQ of odd order n. Then, the set of all l-spins of $Q$ is a cyclic group of odd order n
under composition of mapping denoted as $(\Phi_{L},\circ)$.
\end{myth}
{\bf Proof:}\\
Let $$ P_{R}=\{ \lambda_{i}, \varphi_{j} \in M(Q,\cdot)|\lambda_{ij}=\lambda_{i}\lambda^{-1}_{j}=\lambda_{i}\varphi_{j}, i,j\in Q \} .$$  The remaining part is similar to the proof of Theorem \ref{quad04}.
\begin{myth}
Let $Q(\cdot)$ be an IPQ of order $3^{n}, n\ge 1$. Then, the left and right spins coincide.
\end{myth}
{\bf Proof}\\
By Theorem \ref{quad21} $\lambda_{i}=\varphi_{i}$ and $\lambda_{j}=\varphi_{j}$.\\
Then, consider: $$ \lambda_{ij}=\lambda_{i}\varphi_{j}= \varphi_{i}\varphi_{j}=\varphi_{i}\lambda_{j}=\varphi_{ij} $$
\begin{myrem}
We speak of spins of IPQ since all spins of IPQ are left and right simultaneously.
\end{myrem}

\section{Conclusion}
The concept of middle translation is a 'track algebra' where each element in an algebra is either tracked from the left (left middle translation) or tracked from the right (right middle translation).
This paper, therefore investigated the consequences of left (right) middle translations on involutory latin quandles as well as their induced representations. These permutations $(\lambda_{i}~\&~\varphi_{i})$ were further applied on cyclic groups of odd orders to produce involutory latin quandles. Nice applications of these permutations also helped in reversing the process under r-spins and l-spins of these quandles to recover the cyclic groups earlier used.


\begin{thebibliography}{99}
\bibitem{qua11} V. D. Belousov, {\it Fundamentals of the theory of quasiqroups and loops}, Nauka, Moska (1967) (Russian). Reine Angew. Math. 160,  (1929) 111-130.
\bibitem{phd41} R. H. Bruck, {\it A survey of binary systems}, Springer-Verlag, Berlin-G\"ottingen-Heidelberg, (1966)  185pp.
\bibitem{qua12} V. D. Belousov, {\it On group associated with a quasigroup}, (Russian). Mat. Issled. 4 (1969). no 3, 21-39.
\bibitem{qua13} J. Denes and A. D. Keedwell {\it Latin Squares and their Applications}, Academia Kiado, Budapest (1974).
\bibitem{qua15} I. I. Deriyenko (Derienko), {\it Necessary Conditions of the Isotopy of finite Quasigroups}, (Russians) Mat. Issled. 120 (1991), 51-63.
\bibitem{qua14} I. I. Deriyenko, {\it On middle translations of finite quasigroups}, Quasigroups and Related Systems 16 (2008), 17-24.
\bibitem{qua16} M. Elhamdadi, {\it Distributivity in Quandles and Quasigroups}, Algebra, Geometry and Mathematical
Physics, Springer Proceedings in Mathematics and Statistics,  85,Springer-Valag Heidelberg, (2014) 325-340.
\bibitem{qua17} A. O. Isere, {\it A quandle of order 2n and the concept of quandles isomorphism}, Journal of Nigerian Mathematical Society, 39(2), (2020) 155-166.
\bibitem{qua18} A. O. Isere, O. A. Elakhe and C. Ugbolo , {\it A Higher Quandle of order 24, and its Inner Automorphisms }, J. Physical \& Applied Sciences,  1, 1(2018) 100-110.
\bibitem{qua19} A. O. Isere, J. O. Adeniran and T.G. Jaiyeola, {\it Latin Quandles and Applications to Cryptography }, Mathematics for Applications, 10(2021), 37-53.
\bibitem{phd20} A. O. Isere, J. O. Adeniran and A. A. A. Agboola, {\it Representations of Finite Osborn Loops}, Journal of Nigerian Mathematical Society, 35(2), (2016) 381-389.
\bibitem{qua21} A. O. Isere and J. Ezurike, {\it On Involutory Latin Quandles}, Submitted to International Journal of Group Theory.
\bibitem{qua01} D. Joyce, {\it A classifying invariant of knots, the Knot Quandle }, J. Pure Appl. Alg. 23, 1(1982) 37-66.
\bibitem{qua02} D. Joyce, {\it Simple Quandles }, J. Alg. 79, (1982) 307-318.
\bibitem{qua07} S. Kamada, H. Tamaru and K. Wada, {\it On classification of Quandles of cycle type }, Tokyo Journal
of Mathematics 39 1(2016) 157-171.
\bibitem{qua23} A. Krapez, {\it A Note On Belousov quasiqroups}, Quasigroups and Related System 15,  (2007) 291-294.
\bibitem{qua22} A. D. Keedwell and V. A. Shcherbacov, {\it On m-inverse loops and quasigroups with a long inverse cycles}, Australasian Journal of combinatorics, 26, (2002) 99-119.
\bibitem{qua24} A. D. Keedwell and V. A. Shcherbacov, {\it Quasigroups with an inverse property and generalized parastrophic identities}, Australasian Journal of combinatorics, 13, (2005) 109-124.
\bibitem{qua09} S. Matveev, {\it Distributive groupoids in knot theory }, (Russian) mat. sb. (N. S) 119 (161) (1982), 1, 78-88, 160. \emph{MR672410}.
\bibitem{qua27a} D. A. Stanovsky, A guide to self-distributive quasigroups, or latin quandles, Quasigroups and Related Systems, 23 (2015), 91-128.
\bibitem{qua27} I. Stuhl and P. Vojtechovsky, Enumeration of involutory latin quandles, Bruck loops and commutative automorphic loops of odd prime power order, Nonassociative Mathematics and its applications, 261-276, Contemp. Math. 721, Amer. Math. Soc, Providence, RI, 2019.
\bibitem{qua28} V. A. Shcherbacov, {\it Elements of quasigroup theory and applications}, Boca Raton, Taylor and Francis (2017)
\bibitem{qua30} V. A. Shcherbacov, {\it Some properties of full associated group of IP-loop}, (Russian), Izvestia AN Mold. SSR-Ser fiz-techn. i mat. nauk 2 (1984),51-52.
\end{thebibliography}
\end{document}